\documentclass[12pt]{amsart}

\usepackage{amsmath,amsfonts,amssymb,amsthm}
\usepackage{mathtools}
\usepackage[numbers,sort&compress]{natbib}
\usepackage{graphicx}
\usepackage[dvipsnames]{xcolor}
\usepackage{comment}
\usepackage[T1]{fontenc}
\usepackage{geometry}
\usepackage[colorlinks=true,
            linkcolor=blue,
            citecolor=blue,
            urlcolor=blue]{hyperref}

\usepackage{palatino}
\usepackage{etoolbox}
\AtBeginEnvironment{abstract}{\small}

\definecolor{DarkOliveGreen}{rgb}{0.08,0.42,0.08}

\usepackage{fancyhdr}
\numberwithin{equation}{section}

\newtheorem{thm}{Theorem}[section]

\newtheorem{lemma}[thm]{Lemma}
\newtheorem{df}[thm]{Definition}
\newtheorem{proposition}[thm]{Proposition}

\newcommand{\pNC}{\mathbb{P}^{\mathrm{NC}}}
\newcommand{\pex}{\mathbb{P}^{\mathrm{Ex}}}

\newcommand{\indic}[1]{\mathbf{1}_{#1}}

\title{Exclusion vs Noncolliding dynamics in Markov processes}

\author{Rohan Shiatis}
\address{Department of Mathematics, King's College London, Strand, London, WC2R 2LS, UK}
\email{rohan.shiatis@kcl.ac.uk}

\begin{document}

\begin{abstract}
Let $(\xi_t)_{t \geq 0}$ be a continuous-time Markov process taking values in a countable state space $E$. We consider $k$ copies of $\xi$ started from distinct positions, and write $X = (X_t)_{t \geq 0}$ for the resulting process of configurations in $E^k$. We compare two forms of interaction between the particles: the first is conditioning the particles to never collide (noncolliding dynamics), and the second is simply switching off any jumps that would cause two particles to occupy the same position (exclusion dynamics). We show under mild conditions on the tail probabilities of the collision time that the dynamics are related by an explicit exponential martingale change of measure of one another in the sense that 
    \begin{equation*}
        \left.\frac{d\pNC}{d\pex}\right|_{\mathcal{F}_t} = \frac{h(X_t)}{h(X_0)}\exp\left(\lambda t - \int_0^t \kappa(X_u)du\right),
    \end{equation*}
where $h(x)$ is a principal eigenfunction and $\kappa(x)$ is a killing rate. This generalises a recent result in the integrable probability literature relating TASEP on the ring to asymmetric walkers conditioned to never collide.

\end{abstract}

\maketitle


\vspace{0.5em} \noindent\textbf{Keywords:} Noncolliding Markov processes, interacting particle systems,
totally asymmetric simple exclusion process, Doob $h$-transforms,
determinantal transition kernels, Karlin--McGregor formula.

\section{Introduction and Definitions}

\subsection{A motivational example}
Before discussing our results in full, we begin with an example, discussing a recent observation made in \cite{johnston} stating that the totally asymmetric exclusion process on the ring can be obtained from a change of measure involving a so-called `traffic integral'. 
To discuss this work, let $\mathbb{Z}_n := \mathbb{Z}/n\mathbb{Z} := \{0,1,\ldots,n-1\}$ denote the integers mod $n$. Consider the straightforward Markov process $(X_t)_{t \geq 0}$ on $\mathbb{Z}_n$ defined by the simple dynamics that we jump from $j$ in $\mathbb{Z}_n$ to $j+1$ (mod $n$) at rate $1$, and with no other dynamics. For $1 \leq k \leq n-1$, consider now several different Markov processes on $\mathbb{Z}_n^k$ which involve $k$ copies of the original process undergoing some interaction:
\begin{itemize}
\item \textbf{The process with independent dynamics.} We let the $k$ walkers evolve independently.
\item \textbf{The process with noncolliding dynamics.} We condition on the event that the walkers never collide (in the sense of Doob).
\item \textbf{The process with exclusion dynamics.} We simply ignore any event that would cause two particles to occupy the same position. More specifically, if at some moment the $i^{\text{th}}$ particle is at state $j$, then this particle jumps from $j$ to $j+1$ at rate $1$ if and only if $j+1$ is unoccupied.
\end{itemize}

The two latter processes above have similar dynamics in the sense that if we start from $k$ distinct positions $x = (x_1,\ldots,x_k)$, they both govern processes on $\mathbb{Z}_n^k$ with $k$ stochastically moving particles, no two of which ever occupy the same position. However, their dynamics are slightly different: conditioning the particles to never collide is different from exclusion dynamics. Qualitatively, in the exclusion case, particles are more likely to spend time in the vicinity of one another. 

It is possible to say more, however. The exact relationship between the two processes can be characterised by an exponential martingale change of measure. For $x \in \mathbb{Z}_n^k$ with distinct coordinates, let $\mathbb P_x^{\mathrm{NC}}$ and $\mathbb P_x^{\mathrm{Ex}}$ denote respectively the laws of the processes with noncolliding dynamics and exclusion dynamics started with a configuration of particles at $x = (x_1,\ldots,x_k)$. We have the following result \cite[Theorem 1.12]{johnston}.

\begin{thm}[Johnston \cite{johnston}] \label{thm:j}
Define the \textbf{traffic} of an element $(x_1,\ldots,x_k)$ of $\mathbb{Z}_n^k$ with distinct coordinates to be the number of particles waiting behind a neighbour:
\begin{align*}
\mathrm{Traffic}(x) := \sum_{1 \leq i \neq j \leq k} \indic{ x_{j}- x_i = 1 \bmod n}.
\end{align*}
Then we have
\begin{align*}
\frac{d\mathbb P_x^{\mathrm{NC}}}{d\mathbb P_x^{\mathrm{Ex}}}\Bigg|_{\mathcal{F}_t} = \frac{h(X_t)}{h(X_0)} \exp \left\{ \int_0^t c_{n,k} -\mathrm{Traffic}(X_s) \,\mathrm{d}s \right\},
\end{align*}
where $(\mathcal{F}_t)_{t \geq 0}$ is the underlying filtration for the process and $c_{n,k} = k - \sin(\pi k /n)/\sin(\pi/n)$, and where
\begin{align} \label{eq:vdm}
h(x) := \prod_{1 \leq i < j \leq k} |e^{2 \pi \iota x_j/n} - e^{ 2 \pi \iota x_i/n} |.
\end{align} 
\end{thm}

This result was subsequently generalised by the author and Johnston \cite{JS} to prove a relationship between the asymmetric exclusion process and asymmetric random walkers on $\mathbb{Z}_n$ conditioned to never collide. 

The value of this result lies in the fact that the exclusion process, which is often seen as a desirable object to understand in integrable probability (not least because of convergence to the KPZ fixed point in various scaling limits \cite{corwin, quastelsarkar}), can be obtained from the noncolliding process, which has been shown to have an explicit integrable structure and nice correlation functions.

\medskip
The purpose of the present article is to show that the relationship in Theorem \ref{thm:j} holds far more generally. Indeed, under some fairly mild conditions, we show that Markov processes with exclusion dynamics can be obtained from Markov processes with noncolliding dynamics by a change of measure. 

\color{black}
\subsection{Overview and main result}
    
Let $E$ be a countable state space and let $\xi=(\xi_t)_{t\geq 0}$ be an irreducible, continuous-time pure-jump Markov process on $E$ with bounded rates and finite out-degree at every state. For $k\geq 2$, consider $k$ particles with configurations in $E^k$, and define the collision set
\begin{equation*}
    \Delta
=
\left\{
(x_1,\ldots,x_k)\in E^k:
x_i=x_j\text{ for some }1\leq i<j\leq k
\right\}.
\end{equation*}

We write
\begin{equation*}
E^k_\circ:=E^k\setminus\Delta
\end{equation*}
for the space of collision-free configurations. Further write $\mathbf{P}$ for the law of $k$ independent walkers each moving according to $\xi$, and $\tau$ for the first time that two walkers collide.

Our aim is to compare two different dynamics on $E^k_\circ$. We consider two probability laws:

\begin{enumerate}
\item The law $\pNC$ of $k$ independent copies of $\xi$, though conditioned so that no two particles ever collide.
\item The law $\pex$ of $k$ particles evolving according to the dynamics of $\xi$, except that transitions which would result in a collision are suppressed.
\end{enumerate}

These laws have similar descriptions but exhibit different behaviour. Under $\pNC$, the process is subject to a global constraint which manifests as a repulsion between particles, whereas under $\pex$, the exclusion interaction produces no such repulsion. Analogously to the motivating example, the purpose of this work is to understand the relationship between these two probability laws on the same path space, and to determine under what conditions one may be expressed in terms of the other via a Radon--Nikodym derivative. 

We would like to establish that this change of measure exists under a fairly mild condition on the underlying (single-particle) jump process. Writing $a(T) \sim b(T)$ to denote $\lim_{T \to \infty} b(T)/a(T) = 1$, we introduce the following condition:

\begin{df}\label{df:k-ESS}
We say the process $\xi$ is $k$-separable if there exists $z \in E^k_\circ$ such that for all $x \in E^k_\circ$, and $t \geq 0$, the limit
\begin{align} \label{eq:ax}
\lim_{T \to \infty} \frac{ \mathbf{P}_x( \tau > T-t)}{ \mathbf{P}_{z}(\tau > T) } 
\end{align}
exists. 

\end{df}

Our reason for calling the process $k$-separable is that the following result from Jacka and Roberts \cite[Lemma 2.3]{JR} guarantees that if \eqref{eq:ax} holds, the tail probability of the collision time decomposes as a product of a function of initial position and time.

\begin{lemma}[cf. Jacka and Roberts \cite{JR}]\label{lem:k-ESS}
If \eqref{eq:ax} holds, then there exist functions 
$h:E^k \to [0,\infty)$ with $h(x) = 0$ if and only if $x \in \Delta$, and $a\colon [0,\infty) \to (0,\infty)$ non-increasing satisfying $\lim_{T \to \infty} a(T+t)/a(T) = e^{ - \lambda t}$ for some $\lambda \geq 0$, such that
we have the asymptotics as $T \to \infty$,
\begin{equation} \label{eq:separable}
    \mathbf{P}_x(\tau > T) \sim h(x)a(T) \quad \text{ for } x \in E^k_\circ.
\end{equation}
\end{lemma}

We remark that when $x \in \Delta$, for all $T$, $\mathbf{P}_x(\tau > T) = 0$.
\begin{proof} 
    Set
    \begin{equation*}
        a(T) = \mathbf{P}_{z}(\tau > T),
    \end{equation*}
    a non-increasing function. If we let $b(t) = \lim_{T \to \infty}a(T+t)/a(T)$, then $b(t)$ is a non-increasing function that satisfies
    \begin{equation*}
        b(t)b(s) = b(t+s),
    \end{equation*}
    for all $s,t \geq 0$. Such functional equations have only one class of solutions, namely $b(t) = e^{-\lambda t}$, with $\lambda = -\log b(1) \geq -\log b(0) = 0$.

    From \eqref{eq:ax}, for all $x \in E^k$, we can define
    \begin{equation*}
        h(x) = \lim_{T \to \infty}\frac{\mathbf{P}_x(\tau > T)}{a(T)}.
    \end{equation*}
    If we let
    \begin{equation*}
         E^k_+ = \{x \in E^k : h(x) > 0\} \subset E^k_\circ,
    \end{equation*}
    then this gives \eqref{eq:separable} for $x \in E^k_+$, and we have the required result if $E^k_+ = E^k_\circ$. Note that $h(z) = 1$, so $z \in E^k_+$, and that $h = 0$ on $\Delta$, since $\mathbf{P}_x(\tau > T) = 0$. 
    
    Write that $y \to x$ if $\mathbf{P}_y(X_s = x, \tau > s) > 0$. Supposing now that $y\to x$ and $x \in E^k_+$, then for all $T \geq s$, by the Markov property,
    \begin{equation*}
        \mathbf{P}_y(\tau > T) \geq \mathbf{P}_y(X_s = x, \tau > s)\mathbf{P}_{x}(\tau > T-s).
    \end{equation*}
    Hence, dividing by $a(T)$ and taking $T \to \infty$,
    \begin{equation*}
        h(y) \geq \mathbf{P}_y(X_s = x, \tau > s) h(x)e^{\lambda s} > 0.
    \end{equation*}
    So $y \in E^k_+$ also. In particular, this holds whenever $x = z$.

    Let $\sigma$ be a permutation on $\{1,\dots,k\}$ and for any $x \in E^k$, define
    \begin{equation*}
        x^\sigma = (x_{\sigma(1)}, \dots x_{\sigma(k)}).
    \end{equation*}
    Then observing that for any $x \in E^k_{\circ}$, the law of $\tau$ under $\mathbf{P}_{x^\sigma}$ is the same over all permutations $\sigma$,
    \begin{equation*}
        h(x^\sigma) = \lim_{T \to \infty}\frac{\mathbf{P}_{x^\sigma}(\tau > T)}{a(T)} = \lim_{T \to \infty}\frac{\mathbf{P}_x(\tau > T)}{a(T)} = h(x).
    \end{equation*}
So in particular, it suffices to show that $y \to z^\sigma$ for some permutation $\sigma$, for any $y \in E^k_\circ$. We show this is possible by construction.

Write $z = (z_1, \dots, z_k)$ and suppose, relabelling the entries of $z$ if necessary, that exactly $m < k$ of the entries in $z$, namely $z_1, \dots, z_m$, are occupied by walkers in $y$. Choose a walker $y_r$ not at any state that is an entry of $z$. Since $\xi$ is irreducible, there is a finite directed path $P$ from $y_r$ to $z_{m+1}$, visiting no state twice, with every move along it having positive rate. Say that $P$ passes through the walkers $y_{a_1}, \dots, y_{a_l}$ in this order, and write $p_0 = y_r$, $p_n = y_{a_n}$ for $1 \leq n \leq l$, and $p_{l+1} = z_{m+1}$. Starting with the walker at $p_l$ and working backwards along the directed path, we move the walker at $p_n$ along $P$ to $p_{n+1}$, for $0 \leq n \leq l$. Each such move is collision-free, and the killed process performs this finite sequence of moves with positive probability.

Upon doing so, we realise that the locations with walkers are the same, except for having no walker at $y_r$ and having a walker at $z_{m+1}$, thus increasing the number of walkers in locations in $z$ by 1 (although the walkers occupying $z_1, \dots, z_m$ may have changed). By induction on $m$, we may move all walkers to locations in $z$, so that $y \to z^\sigma$ for some $\sigma$.
\end{proof}

In the setting of Definition \ref{df:k-ESS}, we will see in the sequel that the function $h$ is an eigenfunction of the generator of the process with $k$ independent particles, killed when it hits $\Delta$. In some settings, the function $a$ may be taken to be equal to a constant. We emphasise in particular that the asymptotic decay rate $\lambda$ in Lemma~\ref{lem:k-ESS} may be equal to zero.

Many common processes have this property: for starters, in the setting of Theorem \ref{thm:j} above, we can take $h$ as in \eqref{eq:vdm}, and $a(t) = e^{ - c_{n,k}t}$. 
Generally, if $E$ is finite and irreducible, $a(t)$ will decay exponentially and hence $\lambda$ will be strictly positive.

We obtain different behaviour, however, in infinite domains without strong recurrence, and in these cases $\lambda$ may be zero and $a(t)$ will have polynomial decay or even be constant. 
Take, for instance, $k$ independent continuous-time simple random walks on $\mathbb{Z}^d$. We find that \eqref{eq:separable} holds, with the following regimes according to whether $d = 1, d= 2$, or $d \geq 3$:
\begin{itemize}
\item If $d = 1$, we can take $h(x) = \prod_{1 \leq i < j \leq k} |x_j - x_i|$ and $a(t) = c t^{ - k(k-1)/4}$.
\item If $d = 2$, $h(x)$ is a function that is difficult to describe explicitly, but we have $a(t) = \log(t)^{-\binom{k}{2}}$ \cite{CK} (not rigorously established for $k \geq 3$). 
\item If $d \geq 3$, the underlying random walk is transient (and thus so is the difference between any two walkers). It follows that $h(x) := \mathbf{P}_x(\tau = \infty) > 0$, and we can set $a(t) = 1$ in \eqref{eq:separable}.  
\end{itemize}

In any case, these examples of simple random walks on discrete tori or on the $d$-dimensional lattice are the prototypical examples to keep in mind.

We are now ready to come to our main result:

\begin{thm}\label{thm:main}
Suppose that $\xi$ is $k$-separable (as in Definition \ref{df:k-ESS}). Then, for all $x \in E^k_\circ$, there exists a function $\kappa:E^k_\circ \to \mathbb{R}$ such that $\pNC_x$ and $\pex_x$ are related via the exponential martingale change of measure,
    \begin{equation*}
        \left.\frac{d\pNC_x}{d\pex_x}\right|_{\mathcal{F}_t} = \frac{h(X_t)}{h(X_0)}\exp\left(\lambda t - \int_0^t \kappa(X_u)du\right).
    \end{equation*}

    The function $\kappa(x)$ is the total rate of jumps into $\Delta$ from $x$:
    \begin{align*}
    \kappa(x) := \sum_{y \in \Delta} L(x,y),
    \end{align*}
    where $L$ is the infinitesimal generator for the process with law $\mathbf{P}$.
\end{thm}

Here, both laws are considered as measures on the path space $D([0,\infty),E^k)$. Since $h > 0$ on $E^k_\circ$ by Lemma~\ref{lem:k-ESS}, the density above is strictly positive, so $\pNC_x$ and $\pex_x$ are mutually absolutely continuous on $\mathcal F_t$.

\subsection{Examples}

\subsubsection{TASEP on the ring}

Consider the totally asymmetric exclusion process on the ring, where $E = \mathbb{Z}_n \coloneqq \mathbb{Z}/n\mathbb{Z}$ and $k < n$ walkers each jump $x \to x+1$ at rate 1, but jumps to occupied sites are suppressed. Since $E$ is finite, $\xi$ is $k$-separable.

Since the states are cyclically ordered, the killed process on $E^k_\circ$ is not irreducible, splitting into $(k-1)!$ finite, irreducible and isomorphic components. For each component, $\lambda$ is the same, and $h$ is preserved by the isomorphism. In particular, we find that
\begin{equation*}
    h(x) = \prod_{1 \leq i < j \leq k}\left|e^{2\pi \iota x_j/n} - e^{2\pi \iota x_i/n}\right|, \qquad \lambda = k - \frac{\sin(\pi k/n)}{\sin(\pi/n)},
\end{equation*}
and applying Theorem~\ref{thm:main} recovers \cite[Theorem 1.12]{johnston} as stated in Theorem~\ref{thm:j}.

\subsubsection{ASEP on $\mathbb{Z}^1$}

Consider the asymmetric exclusion process on $\mathbb{Z}$,  where $k < \infty$ walkers each jump $x \to x+1$ at rate $p$ and $x \to x-1$ at rate $q$, with $p,q > 0$, with jumps to occupied sites being suppressed. The rates are bounded, every vertex has out-degree $2$, and one can show that $\xi$ is $k$-separable by analysing the tail asymptotics.

Similar to the earlier discussion of simple random walks on $\mathbb{Z}^d$, we find that
\begin{equation*}
    h(x) = \prod_{1 \leq i < j \leq k}\left|x_j - x_i\right|, \qquad \lambda = 0.
\end{equation*}
Furthermore, we have
\begin{equation*}
    \kappa(x) = (p+q)\mathrm{Traffic}(x),
\end{equation*}
where the traffic is analogously defined as
\begin{equation*}
    \mathrm{Traffic}(x) := \sum_{1 \leq i \neq j \leq k} \indic{ x_{j}- x_i = 1}.
\end{equation*}
giving the change of measure as
\begin{equation*}
    \left.\frac{d\pNC}{d\pex}\right|_{\mathcal{F}_t} = \frac{h(X_t)}{h(X_0)}\exp\left(-(p+q) \int_0^t \mathrm{Traffic}(X_u)du\right).
\end{equation*}
One can realise this formula as the limit of the analogous change of measure for ASEP on the ring as $n \to \infty$, representing how this model is the infinite-volume limit of the finite model.

\subsection{Pathological examples}

The precise conditions we placed on the problem were chosen to encompass many typical examples for $\xi$, such as simple random walks. But it is non-trivial to understand what the most general version of the problem is. For example, we do not know of a rigorous construction of an irreducible $\xi$ with bounded rates and finite out-degree for every $x \in E$ that fails to be $k$-separable, although we expect such examples to exist.

The other conditions have more immediate reasons to be necessary. To outline an informal argument as to why finite out-degree is necessary, imagine a configuration $x \in E^k_\circ$ can jump to neighbours $y_1, y_2,\dots$ under $\mathbf{P}_x$, and the event $\{\tau > T\}$ is most likely achieved by first jumping to some $y_n$ with $n$ which grows with $T$ (even though $\mathbf{P}_{y_n}(\tau > T) \sim h(y_n)a(T)$, this asymptotic need not hold uniformly over $n$). In the limit $T \to \infty$, the process leaves $x$ at a finite rate but lands on no particular state - $\pNC$ loses mass. 

More explicitly, we can give an example with finite out-degrees but allowing for unbounded jump rates, where we construct $\xi$ irreducible and non-explosive, and even have $\pNC$ non-explosive, but have $\pex$ be explosive.

Consider the Markov chain $\xi$ with state space $E = \mathbb{Z}^3 \cup \{1,2,\dots\}$. On $\mathbb{Z}^3$, particles move with rate 1 to each of their $6$ neighbours. We also say that $(0,0,0) \to 1$ at rate 1, and $n \to n+1$ or $n \to (0,0,0)$ both at rate $2^n$ for $n \geq 1$.

Take $k=2$. This process has finite out-degree at every state, but has unbounded rates. Furthermore, if we take $x_0 = (z_1,z_2)$, such that $z_1,z_2,$ and $(0,0,0)$ are sufficiently far apart, then using $\mathbf{P}^\xi$ to denote the law of 1 walker on $E$ and $\mathbb{P}^{SRW}$ to denote the law of $2$ simple random walks on $\mathbb{Z}^3$, by considering whether the particles collide in $\mathbb{Z}^3$ or in $\{1,2,\dots\}$,
\begin{equation*}
    \mathbf{P}_{x_0}(\tau < \infty) \leq \mathbf{P}^\xi_{z_1}(\text{hit } (0,0,0)) + \mathbf{P}^\xi_{z_2}(\text{hit } (0,0,0)) + \mathbb{P}^{SRW}_{x_0}(\text{walkers collide in }\mathbb{Z}^3).
\end{equation*}
But by the transience of a random walk on $\mathbb{Z}^3$, all three quantities on the right-hand side can be made arbitrarily small. So $\mathbf{P}_{x_0}(\tau = \infty) > 0$ and $\frac{\mathbf{P}_x(\tau > T-t)}{\mathbf{P}_{x_0}(\tau > T)}$ converges for all $t \geq 0$, so $\xi$ is $2$-separable.

Furthermore, $\xi$ is transient on $E$, so $\pNC$ certainly exists by conditioning on the event $\{\tau = \infty\}$ which has positive probability, and it inherits non-explosiveness from $\mathbf{P}$. However, if we start with particles at $(0,0,0)$ and $1$, the second particle cannot return to the origin whilst the first particle is there, and so increases along $\{1,2,\dots\}$ with holding times $\mathrm{Exp}(2^n)$, independently of the first particle. In particular, there is a positive probability that the particle at the origin takes time longer than $4$ to jump, and independently, a positive probability that the other particle explodes to infinity in time less than $4$. Thus, $\pex_{((0,0,0),1)}$ explodes with positive probability. 

As a final remark, whilst this example contradicts Theorem~\ref{thm:main} in the sense that the law $\pex$ does not exist on the same path space as $\pNC$, by amending the argument of the paper, one can show that
\begin{equation*}
    \left.\frac{d\pNC_x}{d\pex_x}\right|_{\mathcal{F}_t} = \frac{h(X_t)}{h(X_0)}\exp\left(\lambda t - \int_0^t \kappa(X_u)du\right)\indic{\zeta > t},
\end{equation*}
where $\zeta$ denotes the explosion time and $\pex$ now denotes the law of the minimal exclusion chain up to explosion.

\subsection{The Karlin--McGregor formula}

Write $p_t(x,y)$ for the transition probability of a single copy of the underlying Markov process. Let $x = (x_1,\ldots,x_k)$ and $y = (y_1,\ldots,y_k) \in E^k_\circ$. The Karlin--McGregor formula \cite{KM} states that 
\begin{align} \label{eq:km1}
\det_{1 \leq i,j \leq k} p_t(x_i,y_j) = \sum_{ \sigma \in \mathcal{S}_k} \mathrm{sgn}(\sigma) \mathbf P_x \left( X^i_t = y_{\sigma(i)} ~\forall 1 \leq i \leq k, \tau > t \right).
\end{align}

In some natural cases, many of the terms in the sum in \eqref{eq:km1} will be zero. For instance, if the underlying Markov process is totally ordered, i.e., if $X$ is a jump process such that $X^1_t < \ldots < X_t^k$ for all $0 \leq t < \tau$, and we also have that $y_1 < \ldots < y_k$, the only nonzero term on the right-hand side will occur when $\sigma$ is the identity, so that
\begin{align*}
\det_{1 \leq i,j \leq k} p_t(x_i,y_j)  = \mathbf{P}_x ( X^i_t = y_i , \tau > t ) \qquad \text{(for totally ordered processes)}. 
\end{align*}
This yields an explicit formula for the transition probabilities of $\mathbb{P}^{\mathrm{NC}}$: for $x, y \in E^k_\circ$ with their components ordered,
\begin{align*}
\mathbb{P}_x^{\mathrm{NC}} ( X_t = y ) = e^{ \lambda t} \frac{h(y)}{h(x) } \det_{1 \leq i,j \leq k} p_t(x_i,y_j).
\end{align*}

\subsection{Related work}

Processes conditioned not to collide have been studied extensively in the literature. They often admit explicit formulas involving determinants, as a remnant of the Karlin--McGregor formula above. This means that many such processes can be studied within the framework of determinantal point processes; see, for example, \cite{HKPV}. Dyson \cite{dyson} introduced his eponymous Dyson Brownian motion as a model for the evolution of the eigenvalues of a Hermitian matrix-valued Brownian motion. Equivalently, Dyson Brownian motion can be realised as a collection of Brownian motions conditioned never to collide. It has been studied extensively in the intervening years; see, for example, \cite{grabiner,katoritanemura} and the references therein. Random walks conditioned not to collide are also well studied, appearing in, among others, \cite{konigoconnellroch,eichelsbacherkonig,denisovwachtel}.

Exclusion processes have likewise been studied extensively; see the foundational work of Spitzer \cite{spitzer} and the monographs \cite{liggettIPS,liggett}. Explicit determinantal formulas for the transition probabilities of TASEP were found by Schütz \cite{schutz}. More recent formulas, including formulas for periodic TASEP and TASEP on a ring, were obtained by Baik and Liu \cite{baikliu}. Although exclusion processes and processes conditioned never to collide both evolve on spaces of particle configurations with no multiple occupancy, their dynamics are generally different. In an exclusion process, attempted jumps onto occupied sites are suppressed, whereas conditioning independent particles never to collide produces a Doob transform whose transition rates depend on the corresponding positive harmonic function. The determinantal formulas arising in the two settings nevertheless bear a striking resemblance; see, for example, \cite{diekerwarren}.

It would be interesting to relate our programme to the phenomenon of intertwining. In several important examples, a process of independent particles, its noncolliding Doob transform, and an exclusion- or pushing-type process arise as different projections of a common dynamics on an interlacing array. One particularly illuminating place to see this phenomenon is the work of Warren \cite{warren}, where Dyson Brownian motion is related by an intertwining to a system of reflected Brownian motions on a Gelfand--Tsetlin pattern. Discrete analogues, involving blocking and pushing particle systems and noncolliding random walks, appear in \cite{warrenwindridge,oconnell,borodingorin}.

\medskip
That completes the introduction. 
In the next section we prove our main result.

\color{black}

\section{\texorpdfstring{Proof of Main Result}{Proof of Main Result}}

\subsection{\texorpdfstring{Explicit Law of $\pNC$}{Explicit Law of pNC}}
It is not clear a priori that the law $\pNC$ exists, so the goal of this section is to prove the more specific statement of Proposition~\ref{prop:isprobmeasure} below, which we now work towards stating. Throughout this section, when we state that $\xi$ is $k$-separable, the functions $h(x)$ and $a(t)$, and the value $\lambda$ are those given by Lemma~\ref{lem:k-ESS}. 

Reiterating earlier definitions, let $E$ be a countable state space and let $\xi=(\xi_t)_{t\geq0}$ be an irreducible, continuous-time, pure-jump Markov process on $E$ with infinitesimal generator $L^\xi:E \times E \to \mathbb{R}$. We will further suppose that 
\begin{equation*}
    N^\xi(x) = \{y \neq x: L^\xi(x,y) > 0\}
\end{equation*}
is a finite set for all $x \in E$, and that
\begin{equation*}
    \hat{q} \coloneqq \sup_{x \in E} \sum_{y\in N^\xi(x)}L^\xi(x,y) < \infty.
\end{equation*}

For $k\geq2$, let
\begin{equation*}
\Omega = D([0,\infty),E^k)
\end{equation*}
be the path space of c\`adl\`ag paths taking values in $E^k$. Let
\begin{equation*}
X=(X_t)_{t\geq0},\qquad X_t(\omega)=\omega(t),
\end{equation*}
for $\omega \in \Omega$ be the canonical coordinate process on $\Omega$, and let
\begin{equation*}
\mathcal F_t=\sigma(X_s:0\leq s\leq t),
\qquad
\mathcal A = \bigcup_{t\geq0}\mathcal F_t,
\qquad
\mathcal F_\infty=\sigma(\mathcal A).
\end{equation*}
For $x=(x_1,\ldots,x_k)\in E^k$, let $\mathbf{P}_x$ denote the law on $(\Omega,\mathcal F_\infty)$ under which $X$ consists of $k$ independent copies of $\xi$, with initial distribution $\delta_x$. Let the transition semigroup of $X$ be given by $P_t(x,y)$, for $x,y \in E^k$.

Define the collision set
\begin{equation*}
    \Delta
=
\left\{
(x_1,\ldots,x_k)\in E^k:
x_i=x_j\text{ for some }1\leq i<j\leq k
\right\},
\end{equation*}
and define the first collision time $\tau$ by
\begin{equation*}
    \tau = \inf\{t \geq 0: X_t \in \Delta\}.
\end{equation*}

Let us also define
\begin{equation*}
    E^k_\circ \coloneqq E^k \setminus \Delta.
\end{equation*}
Then we would like to consider the set functions $\mathbf{P}_x(\cdot | \tau = \infty)$ for all $x \in E^k_\circ$ on $\bigcup_{t \geq 0} \mathcal{F}_t$, defined by
\begin{equation*}
    \mathbf{P}_x(A | \tau = \infty) = \lim_{T \to \infty} \mathbf{P}_x(A|\tau > T).
\end{equation*}
The set function above is a natural candidate for $\pNC_x$, but has the issue that it is not defined appropriately on $\mathcal{F}_\infty$. For instance, if $\mathbf{P}_x(\tau = \infty) = 0$, then 
\begin{align*}
    \mathbf{P}_x(\tau = \infty | \tau = \infty) &= \lim_{T \to \infty} \mathbf{P}_x(\tau = \infty|\tau > T)\\
    &= 0,
\end{align*}
but we want $\pNC_x(\tau = \infty) = 1$. Therefore, we will show the following.
\begin{proposition}\label{prop:isprobmeasure}
    Suppose that $\xi$ is $k$-separable, and fix $x \in E^k_\circ$. Then for every $A \in \mathcal{A}$, the limit
    \begin{equation*}
        \mathbf{P}_x(A | \tau = \infty) \coloneqq \lim_{T \to \infty}\mathbf{P}_x(A | \tau > T)
    \end{equation*} 
    exists, and the resulting set function admits a unique extension to a probability measure $\pNC_x$ on $(\Omega,\mathcal{F}_\infty)$, which is supported on $D([0,\infty), E^k_\circ)$. Under $\pNC_x$, the process $X$ is Markov on $E^k_\circ$, with transition probabilities given explicitly by
    \begin{equation*}
        \pNC_x(X_t = y) = e^{\lambda t}\frac{h(y)}{h(x)}\mathbf{P}_x(X_t = y, \tau >t).
    \end{equation*}
\end{proposition}

To prove this proposition, we first show that $k$-separability yields an eigenfunction of the killed generator. Write $L$ for the infinitesimal generator of $X$ under $\mathbf{P}$. For all $x \in E^k_\circ$, let
\begin{equation*}
    L^\Delta f(x) \coloneqq \sum_{y \in E^k_\circ} L(x,y)f(y).
\end{equation*}

Let us further define $q(x) \coloneqq -L(x,x) \leq k \hat{q}$ and 
\begin{equation*}
    N(x) = \{y \neq x : L(x,y) > 0\}.
\end{equation*}
Since each $N^\xi$ is finite, $N(x)$ must also be finite.

\begin{lemma}\label{lem:eigenvecL}
    If $\xi$ is $k$-separable, then $L^\Delta h =  - \lambda h$ on $E^k_\circ$.
\end{lemma}
\begin{proof}
    Fix $x \in E^k_\circ$ and let $T_1$ be the first jump time of $X$. Consider the event $\{\tau > T + t\}$ and condition on whether $T_1 > t$. Using the (Strong) Markov property of $\mathbf{P}_x$,
    \begin{equation*}
        \mathbf{P}_x(\tau > T + t) = e^{-q(x)t}\mathbf{P}_x(\tau > T) + \int_0^t q(x)e^{-q(x)s}\sum_{y \in N(x)}\frac{L(x,y)}{q(x)}\mathbf{P}_y(\tau > T + t - s) ds,
    \end{equation*}
    where the first term is the probability that $\tau > T + t$ and $T_1 > t$, and the second term is an integral over the density of $\tau > T + t$ and $T_1 = s$. Dividing by $a(T)$ and taking $T \to \infty$, we claim that
    \begin{equation}\label{eq:claim}
        e^{-\lambda t} h(x) = e^{-q(x)t}h(x) + \int_0^t e^{-q(x) s} \sum_{y \in N(x)} L(x,y) h(y) e^{-\lambda(t-s)} ds.
    \end{equation}
    Indeed, the left-hand side and the first term of the right-hand side follow immediately from Lemma~\ref{lem:k-ESS}. The second term on the right-hand side features the pointwise limit of the integrand, and the claimed limit holds by the dominated convergence theorem: since $N(x)$ is finite and the rates are bounded, Lemma~\ref{lem:k-ESS} implies that, for $T$ sufficiently large, there exists a constant $\delta$ such that 
    \begin{equation*}
        e^{-q(x)s}\sum_{y \in N(x)}L(x,y)\frac{\mathbf{P}_y(\tau > T + t - s)}{a(T)} \leq 1 \cdot |N(x)| \cdot k\hat{q} \cdot \max_{y \in N(x)}(h(y) + \delta),
    \end{equation*}
    which is our dominating function. If $N(x)$ were infinite, such a uniform bound over $y$ might not exist.

    Continuing from \eqref{eq:claim}, if $\lambda = q(x)$, then \eqref{eq:claim} simplifies to $\sum_{y \in N(x)} L(x,y)h(y) = 0 = (q(x) - \lambda)h(x)$. Otherwise, computing the integral gives
    \begin{equation*}
        (e^{-\lambda t} - e^{-q(x) t}) h(x) = e^{-\lambda t}\sum_{y \in N(x)} L(x,y) h(y) \frac{1 - e^{(\lambda - q(x))t}}{\lambda - q(x)},
    \end{equation*}
    which simplifies to $\sum_{y \in N(x)} L(x,y)h(y) = (q(x) - \lambda)h(x)$ as well. Since the terms in the sum are only over the off-diagonal terms, and since $h = 0$ on $\Delta$,
    \begin{equation*}
        L^\Delta h(x) = L h(x) = (q(x) - \lambda)h(x) - q(x)h(x) = -\lambda h(x). \qedhere
    \end{equation*}
\end{proof}

By Lemma~\ref{lem:eigenvecL}, we can define a new Markov process, via the infinitesimal generator
\begin{equation*}
    L^h(x,y) \coloneqq L(x,y) \frac{h(y)}{h(x)}, \quad L^h(x,x) \coloneqq \lambda - q(x) \quad y \neq x, x,y \in E^k_\circ.
\end{equation*}
The previous lemma guarantees that this generator is conservative on $E^k_\circ$, and with bounded rates, so is non-explosive. For each $x \in E^k_\circ$, let $\mathbb{Q}_x$ be the law of this Markov process on $(\Omega,\mathcal{F}_\infty)$ started from $x$.
\begin{lemma}\label{lem:martingale}
    Define the process
    \begin{equation*}
        M_t \coloneqq e^{\lambda t}\frac{h(X_t)}{h(X_0)}\indic{\tau > t}.
    \end{equation*}
    If $\xi$ is $k$-separable and $x \in E^k_\circ$, then for all $t \geq 0$ and $A \in \mathcal F_t$,
    \begin{equation*}
        \mathbb{Q}_x(A) = \mathbf{E}_x[M_t\indic{A}]. 
    \end{equation*}
\end{lemma}

\begin{proof}
Let $N_t$ be the number of jumps that occur in the time interval $[0,t]$, which is finite almost surely under $\mathbf{P}_x$ and $\mathbb{Q}_x$, and let $0 <  T_1 < \dots < T_{N_t} \leq t$ (well-defined almost surely) be the first $N_t$ jump times. Let 
\begin{equation*}
    \Sigma^t_n = \{(t_1,\dots,t_n): 0 < t_1 < \dots < t_n \leq t\}
\end{equation*}
be the $n$-dimensional simplex. Then consider the set $\mathcal{C}$ of events of the form
\begin{equation}\label{eq:classofsets}
    \{N_t = n, (T_1,\dots,T_n) \in B, X_{T_i} = x_i, 0 \leq i \leq n  \}
\end{equation}
over $n \in \mathbb{N}$, $B$ Borel subsets of $\Sigma_n^t$, $x_1,\dots, x_n \in E^k$ with $x_i \neq x_{i+1}$ and $T_0 = 0$, $x_0 = x$ fixed. These events form a $\pi$-system, and generate $\mathcal{F}_t$. 

Let $A$ be such an event. First note that, if any $x_i \notin E^k_\circ$, then both $\mathbb{Q}_x(A)$ and $\mathbf{E}_x[M_t\indic{A}]$ are zero, from Lemma~\ref{lem:k-ESS}. So we may assume that all $x_i \in E^k_\circ$. Thus,
\begin{equation*}
    \mathbb{Q}_x(A) = \int_{B} \prod_{i=0}^{n-1} \left(\frac{L^h(x_i,x_{i+1})}{q(x_i) - \lambda }\cdot (q(x_i) - \lambda) e^{-(q(x_i) - \lambda)(t_{i+1} - t_i)}\right) \times e^{-(q(x_n) - \lambda)(t-t_n)} dt_1\dots dt_n,
\end{equation*}
since $\frac{L^h(x_i,x_{i+1})}{q(x_i) - \lambda }\cdot (q(x_i) - \lambda) e^{-(q(x_i) - \lambda)(t_{i+1} - t_i)}$ is the probability of jumping $x_i \to x_{i+1}$ multiplied by the density for the inter-arrival times, and $e^{-(q(x_n) - \lambda)(t-t_n)}$ represents the probability of no jumps after time $t_n$. Simplifying and telescoping,
\begin{equation*}
    \mathbb{Q}_x(A) = e^{\lambda t} \frac{h(x_n)}{h(x)}\int_{B} \left(\prod_{i=0}^{n-1} {L(x_i,x_{i+1})} \right) e^{-\sum_{i=0}^{n} q(x_i)(t_{i+1} - t_i)} dt_1\dots dt_n,
\end{equation*}
with the convention that $t_{n+1} = t$. But the same method gives the same result for $\mathbf{E}_x[M_t \indic{A}]$, so they agree on $\mathcal{C}$.

Finally, consider the set $\mathcal{D}$ of events in $\mathcal{F}_t$ such that if $A \in \mathcal{D}$, then $\mathbb{Q}_x(A) = \mathbf{E}_x[M_t\indic{A}]$. This set is a $\lambda$-system containing the $\pi$-system $\mathcal{C}$. In particular, note that the events of $\mathcal{C}$ with $B = \Sigma_n^t$ form a disjoint union of $\Omega$, so by countable additivity, $\mathbb{Q}_x(\Omega) = \mathbf{E}_x[M_t] = 1$, and $\Omega \in \mathcal{D}$. Thus, by Dynkin's Lemma, $\mathcal{D}$ contains $\mathcal{F}_t$. 
\end{proof}

We can now prove the existence and uniqueness of $\pNC$, the method being to show that $\pNC$ and $\mathbb{Q}$ are the same measure, and to inherit the required properties of $\pNC$ from $\mathbb{Q}$, like the Markov property.

\begin{proof}[Proof of Proposition~\ref{prop:isprobmeasure}]
    Fix $t \geq 0$. For $A \in \mathcal{F}_t$, and using the Markov property of $\mathbf{P}_x$,
    \begin{equation}\label{eq:prop1}
        \mathbf{P}_x(A|\tau = \infty) 
        = \lim_{T \to \infty} \frac{\mathbf{P}_x(A,\tau > T)}{\mathbf{P}_x(\tau > T)}
        = \lim_{T \to \infty}\mathbf{E}_x\left[\indic{A,\tau > t} \frac{\mathbf{P}_{X_t}(\tau > T-t)}{\mathbf{P}_x(\tau > T)}\right].
    \end{equation}
Using Lemma~\ref{lem:k-ESS}, the random variable in the expectation converges to $\indic{A,\tau > t} e^{\lambda t}\frac{h(X_t)}{h(x)} = M_t \indic{A}$ as $T \to \infty$, so we would like to justify exchanging the limit and expectation. 

To this end, define the random variable $Y^T_t = \indic{\tau > t} \frac{\mathbf{P}_{X_t}(\tau > T-t)}{\mathbf{P}_x(\tau > T)}$, which we note is non-negative, and satisfies $Y^T_t \to M_t$ as $T \to \infty$. Furthermore, by using the same method as in \eqref{eq:prop1} but in reverse,
\begin{equation*}
    \mathbf{E}_x[Y^T_t] = \mathbf{E}_x\left[\indic{\tau > t} \frac{\mathbf{P}_{X_t}(\tau > T-t)}{\mathbf{P}_x(\tau > T)}\right] = \mathbf{P}_x(\tau > T|\tau > T) = 1
\end{equation*}
so each $Y^T_t$ is integrable, with $\mathbf{E}_x|Y^T_t| \to 1 = \mathbf{E}_x|M_t|$. Note that the last equality comes from Lemma~\ref{lem:martingale}. So applying Scheff\'e's lemma, $Y^T_t$ converges to $M_t$ in $L^1(\mathbf{P}_x)$ and 
\begin{equation*}
    \left|\mathbf{E}_x[\indic{A}Y^T_t] - \mathbf{E}_x[\indic{A}M_t]\right| = \left|\mathbf{E}_x[\indic{A}(Y^T_t - M_t)]\right| \leq \mathbf{E}_x|Y^T_t - M_t| \to 0.
\end{equation*}
So indeed, we have
\begin{equation*}
    \mathbf{P}_x(A|\tau = \infty) 
        =\mathbf{E}_x[M_t \indic{A}].
\end{equation*}
So there exists a measure, namely $\pNC \coloneqq\mathbb{Q}_x$, extending $\mathbf{P}_x(\cdot|\tau = \infty)$ and uniqueness follows e.g.~ by Dynkin's Lemma, since $\mathcal{A}$ is a $\pi$-system generating $\mathcal{F}_\infty$. Finally, for the transition probabilities, using Lemma~\ref{lem:martingale},
\begin{equation*}
    \pNC_x(X_t = y) = \mathbf{E}_x[M_t\indic{X_t = y}] = e^{\lambda t}\frac{h(y)}{h(x)}\mathbf{P}_x(X_t = y, \tau > t). \qedhere
\end{equation*}

\end{proof}

We now turn to the proof of Theorem~\ref{thm:main}. The law $\pex$ is defined as the law on $(\Omega, \mathcal{F}_\infty)$ of the Markov process supported on $D([0,\infty), E^k_\circ)$ with infinitesimal generator given by
\begin{equation*}
    L^{\mathrm{Ex}}(x,y) = \begin{cases}
        L(x,y) &\text{if } x \neq y \in E^{k}_\circ,\\
        L(x,x) + \kappa(x) & \text{if } y=x,
    \end{cases}
\end{equation*}
where $\kappa(x) = \sum_{y \in \Delta}L(x,y)$. Let us further write 
\begin{equation*}
    q^{\mathrm{Ex}}(x) \coloneqq -L^{\mathrm{Ex}}(x,x) = q(x) - \kappa(x).
\end{equation*}
From here, since Proposition~\ref{prop:isprobmeasure} relates the transitions of the noncolliding process to those of the killed process, the proof of Theorem~\ref{thm:main} will consist of relating the transitions of the exclusion process $\pex$ to those of the killed process also.

\begin{proof}[Proof of Theorem~\ref{thm:main}]

We shall aim to prove that, for all $x \in E^k_\circ$, $t \geq 0$ and $\mathcal{F}_t$-measurable functions $f: \Omega \to \mathbb{R}_{\geq 0}$,
\begin{equation}\label{eq:changeofmeas1}
    \mathbf{E}_x\left[f \,\indic{\tau > t}\right] = \mathbb{E}^{\mathrm{Ex}}_x\left[f\, e^{-\int_0^t\kappa(X_u)du}\right].
\end{equation}
Once proven, the required result follows by substituting $f = e^{\lambda t} \frac{h(X_t)}{h(X_0)} \indic{A}$ for any $A \in \mathcal{F}_t$ and applying Proposition~\ref{prop:isprobmeasure}.

We proceed using a monotone class argument, similar to the proof of Lemma~\ref{lem:martingale}. Recall from \eqref{eq:classofsets} the $\pi$-system $\mathcal{C}$ of events of the form
\begin{equation*}
    \{N_t = n, (T_1,\dots,T_n) \in B, X_{T_i} = x_i, 0 \leq i \leq n  \}
\end{equation*}
over $n \in \mathbb{N}$, $B$ Borel subsets of $\Sigma_n^t$, $x_1,\dots, x_n \in E^k$ with $x_i \neq x_{i+1}$ and $T_0 = 0$, $x_0 = x$ fixed. The $\sigma$-algebra generated by these events is $\mathcal{F}_t$. 
Again, let $A$ be such an event. If some $x_i \in \Delta$, then both sides of \eqref{eq:changeofmeas1} vanish, since we would have $\tau \leq t$ on $A$, and $\pex$ never visits $\Delta$, so we can assume all $x_i \in E^k_\circ$. Then
\begin{equation*}
    \mathbf{E}_x\left[\indic{A, \tau > t}\right] = \int_{B} \prod_{i=0}^{n-1} \left(\frac{L(x_i,x_{i+1})}{q(x_{i})} \cdot q(x_{i})e^{-q(x_i)(t_{i+1} - t_{i})} \right)\times e^{-q(x_n)(t-t_{n})}  dt_1\dots dt_n,
\end{equation*}
since $\frac{L(x_i,x_{i+1})}{q(x_{i})} \cdot q(x_{i})e^{-q(x_i)(t_{i+1} - t_{i})}$ is the probability of jumping $x_i \to x_{i+1}$ multiplied by the density for the inter-arrival times, and $e^{-q(x_n)(t-t_{n})}$ represents the probability of no jumps after time $t_n$. This simplifies to
\begin{equation*}
    \mathbf{E}_x\left[\indic{A, \tau > t}\right] = \int_{B} \prod_{i=0}^{n-1} L(x_i,x_{i+1}) e^{-\sum_{i=0}^{n}q(x_i)(t_{i+1} - t_{i}) }  dt_1\dots dt_n,
\end{equation*}
with the convention that $t_{n+1} = t$. But analogously,
\begin{align*}
    \mathbb{E}^{\mathrm{Ex}}_x\left[\indic{A}e^{-\int_0^t\kappa(X_u)du}\right] &= \int_{B} \prod_{i=0}^{n-1} L^{\mathrm{Ex}}(x_i,x_{i+1}) e^{-\sum_{i=0}^{n}q^{\mathrm{Ex}}(x_i)(t_{i+1} - t_{i}) - \kappa(x_i)(t_{i+1} - t_i) }  dt_1\dots dt_n\\
    &= \int_{B} \prod_{i=0}^{n-1} L(x_i,x_{i+1}) e^{-\sum_{i=0}^{n}q(x_i)(t_{i+1} - t_{i})}  dt_1\dots dt_n\\
    &= \mathbf{E}_x\left[\indic{A, \tau > t}\right].
\end{align*}
Finally, let
\begin{equation*}
    \mathcal{D} = \left\{A \in \mathcal{F}_t: \mathbf{E}_x\left[\indic{A} \indic{\tau > t}\right] = \mathbb{E}^{\mathrm{Ex}}_x\left[\indic{A} e^{-\int_0^t\kappa(X_u)du}\right]\right\}.
\end{equation*}
This set is a $\lambda$-system containing the $\pi$-system $\mathcal{C}$. In particular again, the events of $\mathcal{C}$ with $B = \Sigma^n_t$ form a disjoint union of $\Omega$, so by countable additivity, $\mathbf{E}_x\left[\indic{\tau > t}\right] = \mathbb{E}^{\mathrm{Ex}}_x\left[ e^{-\int_0^t\kappa(X_u)du}\right]$, and $\Omega \in \mathcal{D}$. Thus, $\mathcal{D} \supset \sigma(\mathcal{C}) = \mathcal F _t$. So \eqref{eq:changeofmeas1} holds for indicators $\indic{A}$ for $A \in \mathcal{F}_t$. Linearity extends \eqref{eq:changeofmeas1} to simple functions and the monotone convergence theorem extends it to non-negative $\mathcal{F}_t$-measurable functions.
\end{proof}

\section*{Acknowledgements}
The author would like to thank Samuel G. G. Johnston for his support.

\end{document}